\documentclass[leqno,12pt]{amsart} 

\usepackage{amssymb}

\newtheorem{thm}{Theorem}[section]
\newtheorem{prop}[thm]{Proposition}

\newtheorem{cor}[thm]{Corollary}

\theoremstyle{definition}
\newtheorem{defn}[thm]{Definition}

\theoremstyle{remark}

\numberwithin{equation}{section}

\def\C{\mathbb{C}}
\def\Z{\mathbb{Z}}
\def\R{\mathbb{R}}

\def\Z{\mathbb{Z}}
\def\E{\mathcal{E}}
\def\P{\mathbb{P}}

\def\D{\mathcal{D}}

\newcommand{\de}{\partial}
\newcommand{\db}{\overline{\partial}}
\newcommand{\ddb}{{\partial }\overline{\partial}}
\newcommand{\pkk}{\lq\lq$p-$K\"ahler\rq\rq}
\newcommand{\kk}{\lq\lq$1-$K\"ahler\rq\rq}
\newcommand{\nkk}{\lq\lq$(n-1)-$K\"ahler\rq\rq}
\newcommand{\mkk}{\lq\lq$(m-1)-$K\"ahler\rq\rq}
\newcommand{\skk}{\lq\lq$s-$K\"ahler\rq\rq}

\begin{document}






\title[Product of $p-$K\"ahler  manifolds]{Product of generalized  $p-$K\"ahler  manifolds}

\author{Lucia Alessandrini}
\address{ Dipartimento di Matematica e Informatica\newline
Universit\`a degli Studi di Parma\newline
Parco Area delle Scienze 53/A\newline
I-43124 Parma
 Italy} \email{lucia.alessandrini@unipr.it}




\begin{abstract}
A product of K\"ahler manifolds also carries a K\"ahler metric.  In  this short note we would like to study  the  product of generalized $p-$K\"ahler manifolds, compact or not. The results we get extend the known results (balanced, SKT, sG manifolds), and are optimal in the compact case. Hence we can give new non-trivial examples of generalized $p-$K\"ahler manifolds.
\end{abstract}

\maketitle

\section{Introduction}

It is well known  that a product of K\"ahler manifolds  (compact or not)  is a K\"ahler manifold too; as a matter of fact, starting from $(M_1, \omega_1), (M_2, \omega_2)$ with $d \omega_j=0$, we get on $M_1 \times M_2$ the closed strictly positive $(1,1)-$form $\omega := \pi_1^* \omega_1 + 
\pi_2^* \omega_2$, where $\pi_j : M_1 \times M_2 \to M_j$, $j=1,2$, are the standard projections.

The aim of this note is to investigate the product of generalized $p-$K\"ahler  manifolds (see section 3 and \cite{A1}). Some cases are known: in particular that of balanced manifolds (i.e. $(n-1)-$K\"ahler manifolds, where $n$ is the dimension of the manifold)  was considered by Michelsohn (\cite{Mi}), who proved that the product of a couple of balanced manifolds still supports a balanced metric. 

Nevertheless, this result is false for a generic $p$ ($1 < p < n-1$): indeed the situation is very different for $n-$dimensional $p-$K\"ahler manifolds, mainly since we don't speak about a metric, but about a transverse $(p,p)-$form. In particular, let us remark that, when $\omega$ is the K\"ahler form of a K\"ahler metric, then for every $p$, $\omega^p:= \omega \wedge \dots \wedge \omega$ is transverse and closed, but this is no more true when $\omega$ is $\ddb-$closed. And  when $\Omega = \Omega^{p,p}$ is  a $\ddb-$closed and transverse $(p,p)-$form, may be that $\Omega \wedge \Omega$ is not 
$\ddb-$closed nor transverse (see sections 2 and 3).

Our results are collected in section 5: Propositions 5.1 and 5.2 extend the above results to {\it generalized} K\"ahler and balanced manifolds
(as f.i. hermitian symplectic, SKT, strongly Gauduchon manifolds). Next, in Theorems 5.4 and 5.5 we get our main results; moreover, 
we prove by suitable examples that the results are optimal in the compact case. 

Since we need to distinguish among different types of positivity for  vectors, forms, currents, we handle this topic in section 2. In sections 3 and 4 we recall the definition of generalized $p-$K\"ahler manifolds, their characterization by positive currents and some preliminary results. Finally, in section 6, we shall exhibit some classes of examples.
\bigskip

 \section{Positivity}
 
 Let $X$ be a complex manifold of dimension $n \geq 2$, let $p$ be an integer, $1 \leq p \leq n$. 
The aim of this section is to discuss positivity of $(p,p)-$forms, $(p,p)-$vectors and $(p,p)-$currents: we refer to \cite{HK} and to \cite{De} as regards notation and terminology. 
 \smallskip
 
 Positivity  is a local notion, since  it involves only multi-linear algebra; therefore, let us start from a complex $n-$dimensional (euclidean) vector space $E$, its associated (euclidean) vector spaces of $(p,q)-$forms  $\Lambda^{p,q}(E^*)$, 
 and a (orthonormal) basis $\{\varphi_1, \dots, \varphi_n \}$ of $E^*$. (The euclidean structure is not necessary for our purposes, but it is important to assure unicity in some cases, see \cite{HK}).
  
 Let us  denote by $\varphi_I$ the product $\varphi_{i_1} \wedge \dots \wedge \varphi_{i_p}$, where $I = (i_1, \dots, i_p)$ is an increasing multi-index.
 Call $\sigma_p := i^{p^2} 2^{-p}$; thus,  
if $\zeta, \eta \in \Lambda^{p,0} (E^*)$, then $\overline{\sigma_p \zeta \wedge \bar{\eta}} = \sigma_p \eta \wedge \bar{\zeta}$, so that $\sigma_p \eta \wedge \bar{\eta}$ is real;
hence  $\{ \sigma_p \varphi_I \wedge \overline{\varphi_I} , |I| = p \}$ is a (orthonormal) basis of  $\Lambda ^{p,p}_{\R} (E^*):= \{ \psi \in \Lambda ^{p,p} (E^*) / \psi = \overline{\psi} \}$, and 

 $$dv = (\frac{i}{2}  \varphi_1 \wedge \overline{\varphi_1}) \wedge \dots \wedge (\frac{i}{2}  \varphi_n \wedge \overline{\varphi_n}) = \sigma_n \varphi_I \wedge \overline{\varphi_I} , \ I=(1, \dots , n)$$
is a volume form. 
\medskip

We call an $(n,n)-$form $\tau$ {\it positive} ({\it strictly positive}) if $\tau = c\ dv, \ c \geq 0 \ ( c>0)$. We shall write $\tau \geq 0 \ ( \tau > 0)$.

\bigskip

From now on, let $1 \leq p \leq n-1$ and let $k := n-p$.

Remark that every $\Omega \in \Lambda ^{p,p}_{\R} (E^*)$ can be expressed as
 \begin{equation}\label{1}
\Omega = \sigma_p \sum_{|I| = p, |J|=p} \Omega_{I, \overline{J}} \varphi_I \wedge \overline{\varphi_J}
 \end{equation}

where 
$ \Omega_{I, \overline{J}} $ is a hermitian $(N \times N)$ matrix, $N = \frac{n!}{p!(n-p)!}$.

Thus $ \Omega_{I, \overline{J}} $ has real eigenvalues $(\lambda_1, \dots , \lambda_N)$, and there is a (orthonormal) basis $\{\Psi_1, \dots, \Psi_N \}$, 
such that 
 \begin{equation}\label{2}
\Omega = \sigma_p \sum_{j=1}^N \lambda_{j} \Psi_j \wedge \overline{\Psi_j} .
 \end{equation}

This leds to a \lq\lq natural\rq\rq definition of positivity:

\begin{defn}  $\Omega \in \Lambda ^{p,p}_{\R} (E^*)$ is {\it positive} (we shall write: $\Omega \in P^p$) if and only if $\lambda_j \geq 0 \ \forall j$; $\Omega$ is {\it strictly positive} (i.e. $\Omega$ belongs to the interior of $P^p$, $\Omega \in (P^p)^{int}$) if and only if $\lambda_j > 0 \ \forall j$.

\end{defn}

{\it Examples.} The K\"ahler form $\omega$ of a hermitian metric belongs to $(P^1)^{int}$; the $(p,p)-$form $\Omega = \sigma_p \eta \wedge \overline \eta \in P^p$ for every $\eta \in \Lambda^{p,0} (E^*)$.
\medskip

Let us consider the case $p = n-1$: as noticed in \cite{Mi}, p. 279, we have the following result:

\begin{prop} The map $F: (P^1)^{int} \to (P^{n-1})^{int}$ given by $F(\omega) = \frac{1}{(n-1)!}  \omega^{n-1}$ is well-defined and surjective.
\end{prop}

{\it Proof.} If $\omega = \sigma_1 \sum_{j=1}^n \lambda_{j} \psi_j \wedge \overline{\psi_j}$ for some (orthonormal) basis $\{\psi_1, \dots, \psi_n \}$ of $E^*$, then we get easily that
$$\omega^{n-1} =  (n-1)! \ \sigma_{n-1} \sum_{j=1}^n \Lambda_{j} \widehat{\psi_j} \wedge \widehat{\overline{\psi_j}},$$
where $ \widehat{\psi_j} $ means $\psi_J$ with $J = (1, \dots , \widehat{j}, \dots, n)$ and $\Lambda_{j} := \Lambda / \lambda_j$, where $\Lambda := \lambda_1 \cdot \dots \cdot \lambda_n$. 

Thus if all $\lambda_j$ are positive, then all $\Lambda_j$ are positive.

Let $\Omega \in (P^{n-1})^{int}$; then by (2.2) 
$\Omega = \sigma_{n-1} \sum_{j=1}^n \Lambda_{j} \Psi_j \wedge \overline{\Psi_j}$, with 
$\Lambda_{j} > 0$ 
and $\Psi_j =  \widehat{\psi_j} $ for some basis $\{\psi_1, \dots, \psi_N \}$ of $E^*$ (this works because the vector spaces $\Lambda ^{1,1}_{\R} (E^*)$ and $\Lambda ^{n-1,n-1}_{\R} (E^*)$ have the same dimension). 
If we put $$\lambda_j := \frac{(\Lambda_1 \cdot \dots \cdot \Lambda_n)^{\frac{1}{n-1}}}{\Lambda_{j}},$$
and 
$\omega = \sigma_1 \sum_{j=1}^n \lambda_{j} \psi_j \wedge \overline{\psi_j},$ we get by an easy computation that 
$\omega^{n-1} = (n-1)! \ \Omega$.

\bigskip

Let us notice that, in general, whereas in (2.1) the matrix $ \Omega_{I, \overline{J}} $ is not diagonal, in the \lq\lq diagonalized\rq\rq case (2.2) the vectors of  the  basis $\{\Psi_1, \dots, \Psi_N \}$ are not of the type
$ \varphi_{i_1} \wedge \dots \wedge \varphi_{i_p}$, but only linear combinations of elements of this kind. This is the motivation to introduce two more kinds of positivity, as follows:

\begin{defn}   $\eta \in \Lambda^{p,0} (E^*)$ is called {\it simple} (or decomposable) if and only if there are $\psi_1, \dots, \psi_p  \in E^*$ such that $\eta = \psi_{1} \wedge \dots \wedge \psi_{p}$.

$\Omega \in \Lambda ^{p,p}_{\R} (E^*)$ is called {\it strongly positive} ($\Omega \in SP^p$) if and only if 
$ \Omega = \sigma_p \sum_j \eta_j \wedge \overline{\eta_j} ,$ with $\eta_j$ simple.
\end{defn}

\begin{defn}   $\Omega \in \Lambda ^{p,p}_{\R} (E^*)$ is called {\it weakly positive} ($\Omega \in WP^p$) if and only if 
for all $\psi_j \in E^*$, and for all $I = (i_1, \dots, i_k)$ with 
$k+p=n$,
$\Omega \wedge \sigma_k \psi_I \wedge \overline{\psi_I}$ is a positive $(n,n)-$form. It is called {\it transverse} when it is strictly weakly positive, i.e. when $\Omega \wedge \sigma_k \psi_I \wedge \overline{\psi_I}$ is a strictly positive $(n,n)-$form for $\sigma_k \psi_I \wedge \overline{\psi_I} \neq 0$ (i.e. $\psi_{i_1} ,\dots , \psi_{i_k}$ linearly independent).

\end{defn}

As for the link with Definition 2.1, we have the following result:

 \begin{prop} {\rm (see \cite{HK}, Theorem 1.2)} 
$\Omega \in \Lambda ^{p,p}_{\R} (E^*)$ is positive ( $\Omega \in P^p$) if and only if for all $\eta \in \Lambda^{k,0} (E^*)$,
the $(n,n)-$form
$\tau := \ \Omega \wedge \sigma_k \eta \wedge \overline{\eta}$ is positive.
\end{prop}

 {\bf 2.5.1 Remarks.}
 \smallskip
 
 a) Positive forms (in the sense of Definition 2.1) are not considered either by Lelong (\cite{Le}) or by Demailly (\cite{De}); both of them call positive forms (this is the \lq\lq classical sense\rq\rq) what we call weakly positive forms. The strongly positive forms are called {\it decomposable} by Lelong.
 \medskip
 
 b) The sets $P^p, SP^p, WP^p$ and their interior parts are indeed convex cones; moreover, there are obvious inclusions: 
 $ SP^p \subseteq P^p \subseteq WP^p \subseteq \Lambda ^{p,p}_{\R}.$
 \medskip
 
 c) When $p=1$ or $p=n-1$, the three cones coincide, since every $(1,0)-$form is simple (and hence also every $(n-1,0)-$form is simple).
 \medskip 
 
 d) When $1< p< n-1$, the inclusions are strict: indeed, if $\{\varphi_1, \dots, \varphi_4 \}$ is a basis for $\Lambda ^{1,0} (\C^4)$, then it is easy to prove that $\varphi_1 \wedge \varphi_2 + \varphi_3 \wedge \varphi_4$ is not a simple $(2,0)-$form; by Proposition 1.5 in \cite{HK}, this implies that
 $(\varphi_1 \wedge \varphi_2 + \varphi_3 \wedge \varphi_4) \wedge (\overline{\varphi_1 \wedge \varphi_2 + \varphi_3 \wedge \varphi_4})$ is a positive $(2,2)-$form which is not strongly positive.
 
Moreover,  the authors exhibit a $(p,p)-$form which is in the interior of the cone $WP^p$, but has a negative eigenvalue, so it does not belong to the cone $P^p$.
 \medskip
 
 e) Duality. Using the volume form $dv$,  we get the pairing 
 $$f : \Lambda ^{p,p}(E^*) \times \Lambda ^{k,k} (E^*) \to \C$$ 
 given by $f(\Omega, \Psi)dv = \Omega \wedge \Psi$.
 Thus Definition 2.4 can also be stated as:
  $$\Omega \in WP^p \iff \forall \ \Psi \in SP^k, \Omega \wedge \Psi \geq 0.$$
 Moreover, it is not hard to prove that:
 $$\Omega \in SP^p \iff \forall \ \Psi \in WP^k, \Omega \wedge \Psi \geq 0,$$
 $$ \Omega \in P^p \iff \forall \ \Psi \in P^k, \Omega \wedge \Psi \geq 0.$$

\bigskip
As regards vectors, consider $\Lambda_{p,q}(E)$, the space of $(p,q)-$vectors: as before, $V \in \Lambda_{p,0}(E)$ is called a {\it simple vector} if 
$V = v_{1} \wedge \dots \wedge v_{p}$ for some $v_j \in E$; 
in this case, when $V \neq 0$, $\sigma_p^{-1} V \wedge \overline V$ is called a {\it strictly strongly positive} $(p,p)-$vector.

  \begin{prop}
$\Omega \in \Lambda ^{p,p}_{\R} (E^*)$ is  transverse if and only if $\Omega(\sigma_p^{-1} V \wedge \overline V) > 0$
for every $V \in \Lambda_{p,0}(E)$, $V \neq 0$ and simple.
  \end{prop}
  
{\it Proof.} Using the pairing described above, we get an isomorphism $g: \Lambda_{p,p}(E) \to \Lambda ^{k,k} (E^*)$ given as
$f(\Omega, g(A)) = \Omega(A),$ i.e.
$$f(\Omega, g(A)) dv = \Omega \wedge g(A) := \Omega(A) dv, \quad  \forall A \in \Lambda_{p,p}(E), \forall \Omega \in  \Lambda ^{p,p} (E^*).$$
If $\{e_1, \dots, e_n \}$ is a basis of $E$, and $\{\varphi_1, \dots, \varphi_n \}$ is the dual basis, it is easy to check that for all $I = (i_1, \dots, i_p)$, 
$g(\sigma_p^{-1} e_I \wedge \overline {e_I}) = \sigma_k \varphi_J \wedge \overline {\varphi_J}$ with $J =  \{1, \dots, n \} - I.$

Thus 
the isomorphism $g$ transforms $(p,p)-$vectors of the form $\sigma_p^{-1} V \wedge \overline V$, with $V$ simple (i.e. strongly positive vectors), into strongly positive $(k,k)-$forms (of the form $\sigma_k \eta_j \wedge \overline{\eta_j},$ with $\eta_j$ simple). 
Hence we get
$$\Omega(\sigma_p^{-1} V \wedge \overline V) dv = \Omega \wedge g(\sigma_p^{-1} V \wedge \overline V) =  \Omega \wedge \sigma_k \eta \wedge \overline \eta $$
and the statement follows. 

  \medskip

 {\bf 2.6.1 Remark.}  We can identify strictly strongly positive $(p,p)-$vectors  
 (or also forms, using the ismorphism $g$) with $p-$planes in $\C^n$, i.e. with the elements of $G_{\C}(p,n)$; to every plane corresponds a unique unit vector.
 \medskip 
 
 {\bf 2.6.2 Remark.}  It is easy to prove that $\Omega \in WP^p$ if and only if $L^* \Omega$ is a positive form (of maximal degree) for all complex linear maps $L : F \to E$ with ${\rm dim} F = p$ (see \cite{HK} p. 46 and \cite{De}, III.1.6).
  
  \begin{prop} {\rm (see \cite{HK}, Corollary 1.3, \cite{Le}, \cite{De} Ch. III, Proposition 1.12)} 
The  pull-back preserves the different kinds of positivity. That is, when 
  $f$ is a linear map, and  $\Omega \in P^p$ (or $SP^p$, or $WP^p$), then  $f^* \Omega \in P^p$ (or $SP^p$, or $WP^p$).
\end{prop}

 Nevertheless, the wedge product does not preserve {\it weak} positivity:
 
  \begin{prop} {\rm (see \cite{HK}, Corollary 1.3, \cite{Le} Proposition 3, \cite{De} Ch. III, Proposition 1.11)} 
When 
 $\Omega$ and $\Psi$ are positive (or strongly positive), then $\Omega \wedge \Psi$  is positive (or strongly positive); this is no longer true for weakly positive forms, as the cone $SP^p$ is different from the cone $WP^p$ when $1 < p < n-1$. A product of forms, one of which is weakly positive and the other are strongly positive, is weakly positive.
\end{prop}
\medskip

Let us go back to $n-$dimensional manifolds: let $p$ be an integer, $1 \leq p \leq n-1$, and let $k = n-p$; recall that ${\D}^{p,p}(X)_\R$ is the space of compactly supported real $(p,p)-$forms on $X$ and  ${\E}^{p,p}(X)_\R$ is the space of real $(p,p)-$forms on $X$. 

Their dual spaces are: ${\D}_{p,p}'(X)_{\R}$ (also denoted by ${\D '}^{k,k}(X)_{\R}$), the space of real currents of bidimension $(p,p)$ or bidegree $(k,k)$, which we call {\it $(k,k)-$currents}, and 
${\E}_{p,p}'(X)_{\R}$ (also denoted by ${\E '}^{k,k}(X)_{\R}$), the space of compactly supported real $(k,k)-$cur\-rents on $X$. 
Hence a real $(k,k)-$current  is just a real form of bidegree $(k,k)$ with distribution coefficients.

The differential operators $d, \de, \db$ extends naturally to currents by duality; thus we have two De Rham complexes, $(\E^*, d)$ and $((\D')^*, d)$; but the embedding $i: (\E^*, d) \to ((\D')^*, d)$ induces an isomorphism at the cohomology  level. This fact applies also to other cohomologies (as Aeppli and Bott-Chern).

\begin{defn}  
 $\Omega \in {\E}^{p,p}(X)_\R$ is called strongly positive (resp. positive, weakly positive, transverse) if: 
 
 $\forall \ x \in X, \ \Omega_x \in SP^p (T'_xX^*)$ (resp. $P^p (T'_xX^*), \  WP^p (T'_xX^*),$ $(WP^p (T'_xX^*))^{int}$). 

\noindent These spaces of forms are denoted by
$SP^p(X), \ P^p(X),$ $ WP^p(X),$ $(WP^p(X))^{int}$. 
\end{defn}

\begin{defn} ({\rm (see f.i. \cite{HK})}) 
As regards currents in  ${\E}_{p,p}'(X)_{\R}$ we have:

weakly positive currents: $T \in WP_p(X) \iff T(\Omega) \geq 0\ \ \forall \ \Omega \in SP^p(X)$. 

positive currents: $T \in P_p(X) \iff T(\Omega) \geq 0\ \ \forall \ \Omega \in P^p(X)$. 

strongly positive currents: $T \in SP_p(X) \iff T(\Omega) \geq 0\ \ \forall \ \Omega \in WP^p(X)$. 
\end{defn}

 {\bf Notation.} $\Omega \geq 0$ denotes that $\Omega$ is weakly positive;  $\Omega > 0$ denotes that $\Omega$ is transverse; $T \geq 0$ means that $T$ is strongly positive. Thus:
 
 {\bf Claim.} $\Omega > 0$  if and only if $T(\Omega) > 0$ for every $T \geq 0, T \neq 0.$
 \medskip
 
 {\it Examples.} We shall denote by $[Y]$ the current given by the integration on the irreducible analytic subset $Y$ of dimension $p$; $[Y]$ is a closed strongly positive $(n-p,n-p)-$current. 
 By Remark 2.5.1 e), the embedding $i: (\E^*, d) \to ((\D')^*, d)$ maps strongly positive (resp. positive, weakly positive) forms into strongly positive (resp. positive, weakly positive) currents.
 Moreover, let us recall that, if  $f$ is a holomorphic map, and  $T \geq 0$, then  $f_* T \geq 0$.

 \medskip
 
 {\bf Remarks.}  There are obvious inclusions between the previous cones of currents, that is, 
 $SP_p(X) \subseteq  P_p(X) \subseteq WP_p(X)$. The classical positivity for currents (i.e. {\it positive in the sense of Lelong}) is strong positivity; Demailly (\cite{De}, Definition III.1.13) does not consider $P_p(X)$, and indicates $WP_p(X)$ as the cone of positive currents; there is no uniformity of notation in the papers of Alessandrini and Bassanelli.

 \bigskip

\section{ \lq\lq $p-$K\"ahler\rq\rq manifolds}

Let $\Omega \in \E^{n-1,n-1}_{\R} (X)$; by Proposition 2.2, when  $\Omega \in (P^{n-1})^{int}$, then there is $\omega \in (P^1)^{int}$ such that 
 $\Omega = \omega^{n-1}$. 
 Since the proof uses a comparison between the eigenvalues of $\Omega_x$ and those of $\omega_x$, this result is typical of {\it positive} forms. Moreover, when $1< p< n-1$, there are positive forms which are not strongly positive (see Remark 2.5.1 d)): but $\omega^p$ has to be strongly positive by Definition 2.3. Thus this property only holds for $p=n-1$.
 
 Therefore, when $p=n-1$,  it is equivalent to select a strictly weakly positive (i.e. transverse) form, or a strictly positive form, or a strictly strongly positive form, or the $(n-1)-$th power of a K\"ahler form of a hermitian metric. This is not the case when $p < n-1$, so that we will choose in the next Definition  one of the possibilities for forms, i.e. strictly weakly positive forms (see \cite{A1}). 
 \medskip 
 
 To manage in a unified way several different classes of non-K\"ahler manifolds, some of which are well-known in the literature, we introduced in \cite{A1} the following inclusive definition of generalized $p-$K\"ahler manifolds.

 \begin{defn} Let $X$ be a complex manifold of dimension $n \geq 2$, let $p$ be an integer, $1 \leq p \leq n-1$.
 
\begin{enumerate}
\item $X$ is a {\it $p-$K\"ahler (pK) manifold} if it has a closed transverse (i.e. strictly weakly positive) $(p,p)-$form $\Omega \in \E^{p,p}(X)_{\R}$. 

\item $X$ is a {\it weakly $p-$K\"ahler (pWK) manifold} if it has a transverse $(p,p)-$form $\Omega$ with $\de \Omega = \ddb \alpha$ for some form $\alpha$.

\item $X$ is a {\it $p-$symplectic (pS) manifold} if it has a closed transverse  real $2p-$form $\Psi \in \E^{2p}(X)_\R$; that is, $d \Psi = 0$ and $\Omega := \Psi^{p,p}$ (the  $(p,p)-$component of $\Psi$) is transverse.

\item  $X$ is a {\it $p-$pluriclosed  (pPL) manifold} if it has a transverse $(p,p)-$form $\Omega$ with $\ddb \Omega = 0.$
\end{enumerate}

\end{defn}

Notice that:
$pK \Longrightarrow pWK  \Longrightarrow pS  \Longrightarrow pPL;$ 
as regards examples and differences among these classes of manifolds, see \cite{A1}.

When $X$ satisfies one of these definitions, it is called a {\bf generalized $p-$K\"ahler manifold} and in the rest of the paper we will denote it generically as a {\bf \lq\lq$p-$K\"ahler\rq\rq manifold}; obviously, every $n-$dimensional manifold is \lq\lq$n-$K\"ahler\rq\rq, since every $(n,n)-$form is closed.
The form $\Omega$, called a  {\bf \pkk form}, is said to be {\bf \lq\lq closed\rq\rq}. 

\medskip

{\bf 3.1.1 Remark. } As regards Definition 3.1(3), let us write the condition $d \Psi = 0$ in terms of a condition on $\de \Omega$, as in the other statements; 
when $\Psi = \sum_{a+b=2p} \Psi^{a,b},$ then $d \Psi = 0$ is  equivalent to:

i) $\db \Psi^{n-j,2p-n+j} + \de \Psi^{n-j-1,2p-n+j+1}=0$, for $j=0, \dots, n-p-1$, when $n\leq 2p$

and

ii) $\de \Psi^{2p,0}=0, \ \db \Psi^{2p-j,j} + \de \Psi^{2p-j-1,j+1}=0$, for $j=0, \dots, p-1$, when $n > 2p.$

In particular, $\de \Omega = \de \Psi^{p,p}= - \db \Psi^{p+1,p-1}$ (which is the sole condition when $p=n-1$).

On the other hand,  Definition 3.1(2) is equivalent to the existence of a closed transverse real $2p-$form $\Psi = \Psi^{p+1,p-1} + \Psi^{p,p} + \Psi^{p-1,p+1}$, such that $\Psi^{p+1,p-1}=\de \alpha$ for some form $\alpha$. Notice that, for $p=1$,  this is the condition considered in \cite{C2}.

\bigskip

We have just noted that $p=1$ and $p=n-1$ are particular cases; they are more or less known in the literature, with different names: hence we briefly recall them here.

\medskip
{\bf For p = 1}, a transverse form is the fundamental form of a hermitian metric, so that we can speak of $1-$K\"ahler, weakly $1-$K\"ahler,  $1-$symplectic, $1-$pluriclosed {\it metrics}. 

A $1-$K\"ahler manifold is simply a K\"ahler manifold, while
$1-$symplectic manifolds are also called  {\it hermitian symplectic} (\cite{ST}). 
In \cite{Eg} pluriclosed (i.e. $1-$pluriclosed) metrics are defined (see also \cite{ST}); a 1PL metric (manifold) is often called a {\it strong K\"ahler metric (manifold) with torsion} (SKT) (see among others \cite{FPS}). 
Finally, 1WK forms are used, f.i., in \cite{C2}, Theorem 1.2.

\medskip
{\bf For  p = n - 1}, by Proposition 2.2 we get a hermitian metric too. 

This case was studied by Michelsohn  in \cite{Mi}, where $(n-1)-$K\"ahler manifolds are called {\it balanced} manifolds ({\it semi-k\"ahler} in \cite{Ga}). Moreover, $(n-1)-$symplectic manifolds are called {\it strongly Gauduchon manifolds (sG)} by Popovici (compare Definition 3.1(3)  and Theorem 3.2(3) with  \cite{Po1}, Definition 4.1 and Propositions 4.2 and 4.3; see also \cite{Po2}), while $(n-1)-$pluriclosed metrics are called {\it standard} or {\it Gauduchon metrics}. Recently, weakly $(n-1)-$K\"ahler manifolds have been called {\it superstrong Gauduchon (super sG)} (\cite{PU}).

\bigskip

The study of non trivial $p-$K\"ahler manifolds begun with \cite{AA} and \cite{AB1}: 

{\bf when 1 $<$ p $<$ n - 1,} and $\omega$ is a transverse $(1,1)-$form, a standard computation shows that $d \omega^p = 0$ implies $d \omega = 0$ (see also section 6); moreover, a transverse $(p,p)-$form $\Omega$ is not necessarily of the form  $\Omega = \omega^{p}$, where $\omega$ is a transverse $(1,1)-$form, as we explained at the beginning of the section.

Hence in the intermediate cases ($1 < p < n-1$), we must consider a transverse  $(p,p)-$form $\Omega$, which in general  is not of the form $\Omega = \omega^{p}$. Therefore we will not look for \lq\lq good\rq\rq hermitian metrics, but will instead handle transverse forms or strongly positive currents, as we will explain now.  
\medskip

Let us recall the following list of characterization theorems, in connection with Definition 3.1 (see \cite{A1}). They arise from the characterization of K\"ahler ma\-nifolds by currents, which has been introduced by Sullivan \cite{Su} 
and by Harvey and Lawson \cite{HL}. When $T$ satisfies one of the conditions given in the Theorem, we say that $T$ {\bf \lq\lq bounds\rq\rq}. 

\begin{thm} {\bf Characterization of compact generalized $p-$K\"ahler manifolds.}
Let $M$ be a compact complex manifold of dimension $n \geq 2$, and let $p$  be an integer, $1 \leq p \leq n-1$.

\begin{enumerate}
\item $M$ is a {\it $p-$K\"ahler (pK) manifold} if and only if $M$ has no strongly positive currents $T \neq 0$, of bidimension $(p,p)$, such that $T = \de  \overline S + \db S$ for some current $S$ of bidimension $(p,p+1)$ (i.e.  $T$  is the $(p,p)-$component of a boundary). 

\item $M$ is a {\it weakly $p-$K\"ahler (pWK) manifold} if and only if $M$ has no strongly  positive currents $T \neq 0$, of bidimension $(p,p)$, such that $T = \de  \overline S + \db S$ for some current $S$ of bidimension $(p,p+1)$ with $\ddb S = 0$ (i.e.  $T$  is closed and is the $(p,p)-$component of a boundary). 

\item $M$ is a {\it $p-$symplectic (pS) manifold}  if and only if $M$ has no strongly  positive currents $T \neq 0$, of bidimension $(p,p)$, such that $T = d R$ for some current $R$  (i.e.  $T$  is a boundary).

\item  $M$ is a {\it $p-$pluriclosed  (pPL) manifold} if and only if $M$ has no strongly  positive currents $T \neq 0$, of bidimension $(p,p)$, such that $T = \ddb A$ for some current $A$ of bidimension $(p+1,p+1)$. 
\end{enumerate}

\end{thm}

\begin{prop} a) On a compact \pkk manifold $M$, there are no simple exact holomorphic $k-$forms $\alpha \neq 0$.

b) Every compact $n-$dimensional manifold is $(n-1)PL$.
\end{prop}

{\it Proof.} a) If $\alpha$ is such a form, we have $\alpha = \de \beta$, with $\db \beta =0$; thus $\sigma_k \alpha \wedge \overline \alpha \in SP^k(M)$ is $\ddb-$exact, and gives the current $T$ as required in Theorem 3.2(4). Therefore $M$ is not pPL, hence not pK, pWK, pS.

b)  In fact, for $p=n-1$, the current $A$ in Theorem 3.2(4) reduces to a plurisubharmonic global function on a compact complex manifold, hence to a constant. 
 \bigskip

\section{Preliminary results}

\begin{prop} Let $M$ be a \pkk manifold, and $N$ a submanifold of $M$, with ${\rm dim}N > p$;  then $N$ is a \pkk manifold too.
\end{prop}

{\it Proof.} The natural inclusion $i: N \to M$ is a holomorphic immersion, hence it commutes with the operators $d, \de, \db$ and preserves the different kinds of  positivity, by Proposition 2.7. Therefore, when $\Omega$ is a  \pkk form on $M$, then $i^* \Omega$ is a  \pkk form on $N$.

\begin{cor} If $X \times Y$ is a \pkk manifold, and $p < dim X$ (resp. $p < dim Y$),  then $X$ (resp. $Y$) is a \pkk manifold.
\end{cor}

\smallskip
 \begin{prop} Let $M, N$ be  complex manifolds of dimension  $m$ and $n$,  and let $\pi: M \to N$ be a proper holomorphic submersion with $a-$dimensional fibres ($a = m-n$). If $M$ is \pkk for $m > p > a$, then $N$ is \lq\lq $(p-a)-$K\"ahler\rq\rq.
\end{prop}

{\it Proof.} Recall that proper holomorphic submersions allow to push forward forms (see f.i. \cite{De}, I.2.C), so that 
$\pi_* : \E^{p,p}_{\R} (M) \to \E^{p-a,p-a}_{\R} (N)$ and it commutes with the operators $d, \de, \db$.
Moreover, it is not hard to check that  $\pi_*$ preserves the different kinds of  positivity, by definition.
Therefore, when $\Omega$ is a  \pkk form on $M$, then $\pi_* \Omega$ is a  $\lq\lq (p-a, p-a)-$K\"ahler\rq\rq form on $N$.

\medskip

Using the standard projections on factors, we get the following result:

\begin{cor} Let $X \times Y$ be a \pkk manifold; if $Y$ (resp. $X$) is compact and $p > dim Y = n$ (resp. $p > dim X = m$),  then $X$  is a 
\lq\lq $(p-n)-$K\"ahler\rq\rq manifold (resp. $Y$ is a 
\lq\lq $(p-m)-$K\"ahler\rq\rq manifold).
\end{cor}

 \bigskip

\section{Product of  generalized $p-$K\"ahler manifolds}

Let us start recalling a well-known result: {\it The product of K\"ahler manifolds is a  K\"ahler manifold}; this is indeed the principal tool to construct examples of K\"ahler manifolds. The proof of this result also works  for \kk manifolds, i.e. 1WK, 1-symplectic, SKT (1PL) manifolds
(see f.i. \cite{FT2} Remark 2.2).

  \begin{prop} $X, Y$ are \kk manifolds if and only if $X \times Y$ is a \kk  ma\-nifold.
\end{prop}

{\it Proof.} Let $X, Y$ be \kk manifolds,  with \kk forms $\omega_X, \omega_Y$;   consider the standard projections $\pi_X : X \times Y \to X$ and $\pi_Y : X \times Y \to Y$, and take on $X \times Y$ the $(1,1)-$forms  $\pi_X^* \omega_X$ and $\pi_Y^* \omega_Y$. Let $\omega := \pi_X^* \omega_X + \pi_Y^* \omega_Y$.
The  strict positivity of $\omega$ is straightforward as in the K\"ahler case (see also the proof of Theorem 5.4); moreover, $\omega$ is \lq\lq closed\rq\rq since its addenda are \lq\lq closed\rq\rq. The other side is due to Corollary 4.2.

\medskip
The dual case is similar, as we prove in the next Proposition.
 
   \begin{prop} Let $X$ and $Y$ be, respectively, a \mkk manifold of dimension $m$ with \mkk form $\Omega_X$,  and a \nkk manifold of dimension $n$ with \nkk form $\Omega_Y$. Then $X \times Y$ is a \lq\lq$(m+n-1)-$K\"ahler\rq\rq manifold.
   
   Moreover, let $X, Y$ be  complex manifolds of dimension  $m$ and $n$, such that $X$ (resp. $Y$) is compact and the product 
   $X \times Y$ is  \lq\lq$(m+n-1)-$K\"ahler\rq\rq. Then $Y$ is \nkk (resp. $X$ is \lq\lq$(m-1)-$K\"ahler\rq\rq).
\end{prop}

{\it Proof.}   We can suppose $m \geq 3$ or $n \geq 3$ (if not, we are in the \kk case). Recall that there exist hermitian metrics with K\"ahler forms $\omega_X, \omega_Y$ such that $\Omega_X = \omega_X^{m-1}$ and $\Omega_Y = \omega_Y^{n-1}$. As before, 
consider the standard projections $\pi_X : X \times Y \to X$ and $\pi_Y : X \times Y \to Y$, and take on $X \times Y$ the $(1,1)-$forms given by $\pi_X^* \omega_X$ and $\pi_Y^* \omega_Y$. Let $\omega := \pi_X^* \omega_X + \pi_Y^* \omega_Y$.

The  $(1,1)-$form $\omega$ is strictly positive as in the K\"ahler case, so that also $\Omega := \omega^{m+n-1}$ is strictly positive; let us check the  \lq\lq closure\rq\rq. Since
$$\omega^{m+n-1} = c_1 (\pi_X^* \omega_X)^{m-1} \wedge (\pi_Y^* \omega_Y)^n + c_2 (\pi_X^* \omega_X)^{m} \wedge (\pi_Y^* \omega_Y)^{n-1}$$ 
for some positive constants $c_1, c_2$, we get
$$ \de \omega^{m+n-1} = c_1 \de (\pi_X^* \omega_X)^{m-1} \wedge (\pi_Y^* \omega_Y)^n + c_2 (\pi_X^* \omega_X)^{m} \wedge \de (\pi_Y^* \omega_Y)^{n-1},$$ 
$$\ddb \omega^{m+n-1} = c_1 \ddb (\pi_X^* \omega_X)^{m-1} \wedge (\pi_Y^* \omega_Y)^n + c_2 (\pi_X^* \omega_X)^{m} \wedge \ddb (\pi_Y^* \omega_Y)^{n-1}$$ 
and this proves the statement. The second assertion arises from Corollary 4.4.
\medskip

We have just proved that the product of balanced manifolds is balanced (\cite{Mi}, Proposition 1.9)
and moreover that the same holds also for strongly Gauduchon manifolds, Gauduchon manifolds and in the weakly K\"ahler case.
\bigskip

For non trivial $p-$K\"ahler manifolds,  
i.e. when $1 < p < n-1$,  similar simple results are false, also in the compact case and when one of the manifolds is K\"ahler.  
Let us give here   very simple examples of this fact.
\medskip

{\bf 5.2.1 Examples.} Let $Y = I_3,$ the Iwasawa manifold (see section 6), which is not K\"ahler but is 2-K\"ahler (i.e. balanced); 
if $\P_1 \times I_3$ were 2K, 
by Corollary 4.4 we would get that $I_3$ is K\"ahler.
The same holds for  
$(\P_1)^3 \times I_3$: it cannot be 2K.

But, while
$\P_1 \times I_3$ is not 2-K\"ahler, on the contrary $\C \times I_3$ is 2-K\"ahler, as the following Theorem shows (the proof of this result, which is quite technical, can be found in \cite{A2}, Theorem 4.1). 
 
\begin{thm} Let $M$ be a compact holomorphically parallelizable manifold. If $M$ is $p-$K\"ahler, then for every $n \geq 1$, $\C^n \times M$ is $p-$K\"ahler.
\end{thm}

\medskip

In any case,  when $1 < p < dim X-1$, the situation is quite different, because we don't have a metric, but only a transverse $(p,p)-$form. Nevertheless, we can gain some interesting results (we got  a first result  in \cite{AB1}, Theorem 4.5).

\begin{thm} 
Let $X$ be a $m-$dimensional complex manifold which is \skk for all $s, \ p \leq s < m$, and let $Y$ be a $n-$dimensional complex manifold which is \skk for all $s, \ q \leq s < n$. 

Suppose  that $X$ and $Y$ have \skk forms, respectively $\{ \Omega_p, \dots , \Omega_{m-1} \}$ and $\{ \Phi_q,$ 
$\dots , \Phi_{n-1} \}$, with non-negative eigenvalues (i.e. in $(WP^s)^{int} \cap P^s$) (this condition is obvious for $\Omega_{1}, \Omega_{m-1}, \Phi_{1}, \Phi_{n-1}$). 

Let us consider $m-p$ and $n-q$:

If $m-p \leq n-q$, suppose that for every index $a, \ 1 \leq a < m-p$,  either $\de \Omega_{p+a}=0$ or $\de \Phi_{r}=0$ for every $r \in \{n-a, \dots, n-1\}$ (this condition is obvious in the case sK). 

If $m-p > n-q$, suppose that for every index $a, \ 1 \leq a < n-q$,  either $\de \Phi_{q+a}=0$ or $\de \Omega_{r}=0$ for every $r \in \{m-a, \dots, m-1\}$ (this condition is obvious in the case sK). 
\smallskip

Then  $X \times Y$ is \skk for all $s, \ m+n-min(m-p,n-q) \leq s < m+n$, with forms in $(WP^s)^{int} \cap P^s$. 
\end{thm} 

{\bf Remark.} For $p=m-1, q=n-1$ we get exactly Proposition 5.2, while Proposition 5.1 gives more informations.
\medskip

Before to give the proof of this Theorem, we would like to explain the hypotheses regarding positivity on the \skk forms. Recall that we are in the non-trivial \pkk cases, so that we  consider transverse forms, i.e. strictly weakly positive forms. In Theorem 5.4 we ask furthermore that our forms have also non-negative eigenvalues. There is a wide class of  manifolds having those kind of \pkk forms: compact holomorphically parallelizable manifolds (see section 6). 

Moreover we need that the manifolds are \pkk not only for a single index $p$, but from an index $p$ up to the dimension of the manifold: also this very natural property is enjoyed by compact homomorphically parallelizable manifolds (and also by all the known non-trivial examples, see 6.6). One can consider also the class of complex Lie groups with left-invariant $p-$K\"ahler forms (see \cite{A2}).
\medskip

{\it Proof of Theorem 5.4.} It is useful to consider also top-degree forms, i.e. volume forms for $X$ and for $Y$; to use the same notation, pick on $X$ the top degree (closed) form $\Omega_m > 0$, and  on $Y$ the (closed) form
$\Phi_n >0$.

Notice moreover that the hypotheses on non-negative eigenvalues are verified by definition in the cases $p=1, m-1, m$ and $ q = 1, n-1, n$.

Suppose $m-p \leq n-q$ (the other case is similar), and fix an index $j, \ n+p \leq j < n+m$, say $j = n+p+l$ with $0 \leq l < m-p$: we want to prove that  $X \times Y$ is \lq\lq$j-$K\"ahler\rq\rq.

Consider the standard projections $\pi_X : X \times Y \to X$ and $\pi_Y : X \times Y \to Y$, and consider on $X \times Y$ the forms  $\pi_X^* \Omega_s$ and $\pi_Y^* \Phi_s$, which for simplicity are called $\Omega_s$ and $\Phi_s$. 
Let $$\Theta_j := \Omega_{p+l} \wedge \Phi_{n} + \Omega_{p+l+1} \wedge \Phi_{n-1} + \dots + \Omega_{m} \wedge \Phi_{j-m}.$$

{\it First step: Positivity.} Every  addendum belongs to $P^j$, by Proposition 2.8, hence also $\Theta_j$ belongs the convex cone $P^j$.
We will show that $\Theta_j$ is also transverse, using Proposition 2.6, so that $\Theta_j \in (WP^j)^{int} \cap P^j$.

Fix $(x,y) \in X \times Y$, take a simple (unit) $j-$vector $V \in \Lambda_{j,0}(T'_{(x,y)}(X \times Y))$;  $\sigma_j^{-1} V \wedge \overline V$ can be identified (by Remark 2.6.1) with a $j-$plane (also called $V$) in $T'_{(x,y)}(X \times Y)$.
If $a:= dim (V \cap T'_x X)$ and $b:= dim (V \cap T'_y Y)$, we get 
$$j \geq a+b, \ \  p+l = j-n \leq a \leq m, \ \ n-m+p+l = j-m \leq b \leq n.$$
Let us choose now a nice basis for $T'_{(x,y)}(X \times Y) \simeq T'_{x}(X) \oplus T'_{y}(Y)$. Let $(e_1, \dots , e_a)$ be a basis of $V \cap T'_x X,$ and complete it to a basis 
$(e_1, \dots , e_m)$ of  $T'_x X;$ let $(f_1, \dots , f_b)$ be a basis of $V \cap T'_y Y,$ and complete it to a basis 
$(f_1, \dots , f_n)$ of  $T'_y Y.$  

When $j > a+b$, let
$$V = e_1 \wedge \dots \wedge e_a \wedge f_1 \wedge \dots \wedge f_b \wedge v_{a+b+1} \wedge \dots \wedge v_j ,$$
where 
$$v_k := \sum_{i=a+1}^m c_{k,i} e_i + \sum_{i=b+1}^n g_{k,i} f_i.$$ 

{\it Claim.} Eventually rescaling some vectors of the bases, $V$ contains one of the addenda 
$V' := e_1 \wedge \dots e_t \wedge f_1 \wedge \dots f_h,$ with $t+h=j, \ \ p+l \leq t \leq m, \ \ j-m \leq h \leq n.$

{\it Proof of the Claim.} Notice that the indices $(a,b)$ have the same bounds as the indices $(t,h)$; but $t+h=j$, while $a+b \leq j$.

If $a+b=j$, $V' = V = e_1 \wedge \dots \wedge e_a \wedge f_1 \wedge \dots \wedge f_b$ ($t=a, h=b$) satisfies the claim.

When $a+b < j$, so that the vectors $v_k$ really appear in the expression of $V \neq 0$, either one of the $\{e_{a+1}, \dots , e_m\}$ or one of the 
$\{f_{b+1}, \dots , f_n\}$ appear in each $v_k, \ a+b+1 \leq k \leq j$ (that is, at least one of the coefficients $c_{k,i}$ or $g_{k,i}$ does not vanish).
So we can eventually rescale the vectors in the bases $(e_1, \dots , e_m)$ and $(f_1, \dots , f_n)$, to get an addendum $V'$ in $V$ as required.
\medskip

Recall that $\Theta_j = \sum \Omega_{s} \wedge \Phi_{k}$, with $s+k=j, \ \ p+l \leq s \leq m, \ \ j-m \leq k \leq n;$
since every addendum $\Omega_{s} \wedge \Phi_{k}$ of $\Theta_j $ belongs to $ P^j$, we get  $(\Omega_{s} \wedge \Phi_{k}) (\sigma_j^{-1} V \wedge \overline V) \geq 0$; but, by the claim, 
$\sigma_j^{-1} V' \wedge \overline V'$ is an addendum of $\sigma_j^{-1} V \wedge \overline V$, and hence,
when $t=s$ and $h=k$, 
$$(\Omega_{s} \wedge \Phi_{k}) (\sigma_j^{-1} V \wedge \overline V) = (\Omega_{s} \wedge \Phi_{k}) (\sigma_j^{-1} V' \wedge \overline V') > 0;$$

thus $\Theta_j(\sigma_j^{-1} V \wedge \overline V) > 0$.
\medskip

{\it Second step: Closure.} Consider
$$\Theta_j := \Theta_{j_1} + \Theta_{j_2} := (\Omega_{p+l} \wedge \Phi_{n} + \Omega_{m} \wedge \Phi_{j-m}) + 
(\sum \Omega_{s} \wedge \Phi_{k}).$$
where the indices $(s,k)$ in the last sum satisfy: $s+k=j, \  p+l < s < m, \  j-m < k < n.$
We get 
$$\de \Theta_{j_1} = \de \Omega_{p+l} \wedge \Phi_{n} + \Omega_{m} \wedge \de \Phi_{j-m}, \ \ \  \de \Theta_{j_2} =  
\sum (\de \Omega_{s} \wedge \Phi_{k} + \Omega_{s} \wedge \de \Phi_{k}).$$

Thus when both manifolds are $s-$K\"ahler, also $\Theta_j$ is closed. 
\medskip

Let us go to the case sWK: we have 
$\de \Omega_s = \ddb \alpha_s,\  \de \Phi_k = \ddb \beta_k$ for suitable forms $\alpha_s$ and $\beta_k$. In particular, $\ddb \Omega_s = \ddb \Phi_k = 0$. 

We get $\de \Theta_{j_1} = \ddb (\alpha_{p+l} \wedge \Phi_{n} + \Omega_{m} \wedge  \beta_{j-m})$. Moreover, consider a generic summand  
$\Omega_{s} \wedge \Phi_{k}$ of $\Theta_{j_2}$; we get easily
$$\de \Omega_{s} \wedge \Phi_{k} + \Omega_{s} \wedge \de \Phi_{k} = \ddb \alpha_s \wedge \Phi_{k} + \Omega_{s} \wedge \ddb \beta_k = $$
$$ = \ddb (\alpha_s \wedge \Phi_{k} + \Omega_{s} \wedge \beta_k) + (-  \db \alpha_s \wedge \de \Phi_{k} + \de  \alpha_s \wedge \db \Phi_{k} + \db \Omega_{s} \wedge \de \beta_k - \de \Omega_{s} \wedge \db \beta_k).$$
When $\de \Omega_s = 0$, also $\db \Omega_s = 0$ and we can choose $\alpha_s =0$; 
when $\de \Omega_s \neq 0$, by the hypothesis we can take  $\de \Phi_k = 0$ and  $\beta_k =0$. In both cases,  
$$ \de (\Omega_{s} \wedge \Phi_{k} ) = \ddb (\alpha_s \wedge \Phi_{k} + \Omega_{s} \wedge \beta_k).$$
Hence also $\de \Theta_j$ can be expressed as $\ddb \alpha$, for a suitable form $\alpha$ on $X \times Y$.
\medskip

The case sPL is very similar, since as before 
$\ddb \Theta_{j_1} = \ddb \Omega_{p+l} \wedge \Phi_{n} + \Omega_{m} \wedge \ddb \Phi_{j-m}$. Moreover, for a generic summand  
$\Omega_{s} \wedge \Phi_{k}$ of $\Theta_{j_2}$, we get
$\ddb (\Omega_{s} \wedge \Phi_{k}) = - \db \Omega_{s} \wedge \de \Phi_{k} + \de \Omega_{s} \wedge \db \Phi_{k}$, thus 
$\ddb \Theta_j =0$.
\medskip

The case pS requires some computations. We use the characterization of pS manifolds given in Remark 3.1.1, that is, 
our hypothesis on forms $\Omega_s =  \Psi^{s,s}$ and $\Phi_k = \Gamma^{k,k}$ is as follows:

i) $\db \Psi^{m-i,2s-m+i} + \de \Psi^{m-i-1,2s-m+i+1}=0$, for $i=0, \dots, m-s-1$, when $m\leq 2s$

and

ii) $\de \Psi^{2s,0}=0, \ \db \Psi^{2s-i,i} + \de \Psi^{2s-i-1,i+1}=0$, for $i=0, \dots, s-1$, when $m > 2s.$

Moreover, 

i) $\db \Gamma^{n-i,2k-n+i} + \de \Gamma^{n-i-1,2k-n+i+1}=0$, for $i=0, \dots, n-k-1$, when $n\leq 2k$

and

ii) $\de \Gamma^{2k,0}=0, \ \db \Gamma^{2k-i,i} + \de \Gamma^{2k-i-1,i+1}=0$, for $i=0, \dots, k-1$, when $n > 2k.$
\smallskip

Recall that
$$\Theta_j := \Theta_{j_1} + \Theta_{j_2} := (\Omega_{p+l} \wedge \Phi_{n} + \Omega_{m} \wedge \Phi_{j-m}) + 
(\sum \Omega_{s} \wedge \Phi_{k}).$$
where the indices $(s,k)$ in the last sum satisfy: $s+k=j, \  p+l < s < m, \  j-m < k < n.$
As regards $\Theta_{j_1}$ we get 
$$\de \Theta_{j_1} = \de \Omega_{p+l} \wedge \Phi_{n} + \Omega_{m} \wedge \de \Phi_{j-m} = - \db \Psi^{p+l+1,p+l-1}
\wedge \Phi_{n} - \Omega_{m} \wedge  \db \Gamma^{j-m+1,j-m-1} =$$
$$= - \db (\Psi^{p+l+1,p+l-1}
\wedge \Phi_{n} + \Omega_{m} \wedge  \Gamma^{j-m+1,j-m-1}) = - \db A^{j+1,j-1}$$ 
where
$$ - \de A^{j+1,j-1} = - \de (\Psi^{p+l+1,p+l-1}
\wedge \Phi_{n} + \Omega_{m} \wedge  \Gamma^{j-m+1,j-m-1}) = $$
$$ = - \de \Psi^{p+l+1,p+l-1}
\wedge \Phi_{n} - \Omega_{m} \wedge  \de \Gamma^{j-m+1,j-m-1} =$$
$$=  \db \Psi^{p+l+2,p+l-2}
\wedge \Phi_{n} + \Omega_{m} \wedge  \db \Gamma^{j-m+2,j-m-2} =$$
$$ \db (\Psi^{p+l+2,p+l-2} \wedge \Phi_{n} + \Omega_{m} \wedge   \Gamma^{j-m+2,j-m-2}) = \db A^{j+2,j-2}$$
and so on.

As regards $\Theta_{j_2}$ we get
$\de \Theta_{j_2} =  
\sum (\de \Omega_{s} \wedge \Phi_{k} + \Omega_{s} \wedge \de \Phi_{k}).$
Let us consider each summand: by the hypothesis, suppose $\de \Phi_k = 0$. Then:
$$\de (\Omega_{s} \wedge \Phi_{k}) = \de \Omega_{s} \wedge \Phi_{k} = - \db \Psi^{s+1,s-1} \wedge \Phi_{k} =
- \db (\Psi^{s+1,s-1} \wedge \Phi_{k}) = - \db B^{j+1,j-1},$$
where 
$$ - \de B^{j+1,j-1} = - \de (\Psi^{s+1,s-1} \wedge \Phi_{k}) = - \de \Psi^{s+1,s-1} \wedge \Phi_{k} =$$
$$= \db \Psi^{s+2,s-2} \wedge \Phi_{k} = \db (\Psi^{s+2,s-2} \wedge \Phi_{k}) = \db B^{j+2,j-2}$$
and so on. The same holds when $\de \Omega_s = 0$.
Notice that it doesn't matter if we end with a condition of type $\de \Psi^{2s,0}=0$ or not.

Summing up, we get that $\Theta_{j}$ satisfies the conditions of Remark 3.1.1, hence it is a $jS-$form for $X \times Y$.

\bigskip
When $X$ is K\"ahler, we improve Theorem 5.4 in two directions: we don't need the hypothesis on the eigenvalues of the \skk forms on $Y$; moreover, when $m-1 < n-q$, i.e. $m+q < n+1$, we get that $X \times Y$ is \skk also for $m+q \leq s < n+1$.

\begin{thm} 
Let $X$ be a $m-$dimensional K\"ahler manifold, and let $Y$ be a $n-$dimensional complex manifold which is \skk for all $s, \ q \leq s < n$. 

Then $X \times Y$ is \skk for all $s, \ m+q \leq s < m+n$. 
\end{thm} 

{\it Proof.} Choose  \skk forms $\Phi_s$ on $Y$ as before (also $\Phi_n$); recall that on $X$ we can choose
 $\Omega_s := \omega^s$ for a closed transverse $(1,1)-$form $\omega$. Consider their pull-back to $X \times Y$, also called $\Phi_s$ and 
 $\omega$.
 
 Fix an index $j, \ m+q \leq j < m+n$, call  $h = max(0, j-n)$, and let 
 $$\Theta_j := \omega^{h} \wedge \Phi_{j-h} + \omega^{h+1} \wedge \Phi_{j-h-1} +  \dots + \omega^{m} \wedge \Phi_{j-m} =
 \sum_{h \leq k \leq m} \omega^k \wedge \Phi_{j-k}.$$

{\it First step: Positivity.} Every  addendum belongs to $WP^j$, by Proposition 2.8, hence also $\Theta_j$ belongs to the convex cone $WP^j$.
We will show that $\Theta_j$ is also transverse, using Proposition 2.6.

As before, fix $(x,y) \in X \times Y$, take a simple (unit) $j-$vector $V \neq 0$ and identify   $\sigma_j^{-1} V \wedge \overline V$  with a $j-$plane (also called $V$) in $T'_{(x,y)}(X \times Y)$.
If $a:= dim (V \cap T'_x X)$ and $b:= dim (V \cap T'_y Y)$, we get  
$$j \geq a+b, \ \   j-n \leq a \leq m, \ \ j-m \leq b \leq n.$$
Then choose suitable bases 
$(e_1, \dots , e_m)$ of  $T'_x X$ and 
$(f_1, \dots , f_n)$ of  $T'_y Y,$  such that, 
when $j=a+b$, 
$$V = e_1 \wedge \dots \wedge e_a \wedge f_1 \wedge \dots \wedge f_b$$
and when $j > a+b$, for some couple $(k,j-k)$ with $h \leq k \leq m$, a suitable expression of $V$ contains the addendum 
$V' := e_1 \wedge \dots e_k \wedge f_1 \wedge \dots f_{j-k}.$

Since every addendum $\omega^{k} \wedge \Phi_{j-k}$ of $\Theta_j$ belongs to $WP^j$, we get  $(\omega^{k} \wedge \Phi_{j-k}) (\sigma_j^{-1} V \wedge \overline V) \geq 0$; but
$\sigma_j^{-1} V' \wedge \overline V'$ is an addendum of $\sigma_j^{-1} V \wedge \overline V$, and hence 
$$(\omega^{k} \wedge \Phi_{j-k}) (\sigma_j^{-1} V \wedge \overline V) = (\omega^{k} \wedge \Phi_{j-k}) (\sigma_j^{-1} V' \wedge \overline V') > 0;$$
thus $\Theta_j(\sigma_j^{-1} V \wedge \overline V) > 0$.

{\it Second step: Closure.} Since $\de \omega = 0$, by easy computations we get that $\Theta_j$ is \lq\lq closed\rq\rq.
\medskip

\begin{prop}
When $(p,q) \neq (1,1)$, the bounds in the previous Theorems are sharp if and only if both manifolds are compact.
\end{prop}

{\it Proof.} In fact, as we noticed in 5.2.1, $\C \times I_3$ is 2-K\"ahler, while the bound is $m+q =3$ in Theorem 5.5 and $m+n-min(m-p,n-q) = 4$ in Theorem 5.4.
On the contrary, suppose $X$ and $Y$ are compact manifolds. In Theorem 5.4, suppose $m-p \leq n-q$ so that the bound for  $j$ is $n+p \leq j.$ If $X \times Y$ were \lq\lq$(n+p-1)-$K\"ahler\rq\rq, by Corollary 4.4  $X$ would be  \lq\lq$(p-1)-$K\"ahler\rq\rq. The same holds with the bound in Theorem 5.5, using $\pi_Y: X \times Y \to Y$.
 \bigskip

  \section{Examples}
  
Every complex curve is K\"ahler; the situation is quite different in the case of complex surfaces. As regards the compact ones, we have the following picture.

{\bf 6.0 Remark.} Every compact surface is 1PL (SKT), by Proposition 3.3; moreover, there is only a class of special surfaces, because:

1K $\iff$ $b_1$ is even (\cite{La}) 
$\iff$ 1S (\cite{Ga1}, Lemme II.3, or \cite{Po1}, p. 18)  

(in \cite{HL}, Proposition 25, the authors proved that 1S $\iff$ 1WK).

The Hopf surface is not in this class.

Let us notice that this regards manifolds, but not metrics, as it involves the non-existence of currents!

\bigskip

As for dimension bigger than two, first of all we shall  discuss about {\bf positivity}.  
 
 Let us turn back to Definition 3.1: notice that we have chosen the widest cone of forms, that of transverse forms,  because this was the classical definition of strictly {\it positive} forms, and motivated by some geometrical considerations, in particular Remark 2.6.2 and
 what occurs on compact holomorphically parallelizable manifolds: if one of them is \pkk but not K\"ahler, the form $\Omega$ can be chosen in the interior of the cone  $WP^p$, and in $P^p$, but not in the interior of $P^p$ (see the forthcoming Claim 6.4).
 
\medskip
 Let $M$ be a compact holomorphically parallelizable manifold of dimension $n$, that is, $M$ has $n$ holomorphic $1-$forms 
 $\{\varphi_1, \dots, \varphi_n\}$ which are linearly independent at every point; 
  hence $\{\varphi_A\}, |A| = k$, is a basis for $\Omega^k(M)$.

\begin{defn}  If a $p-$K\"ahler (pK) manifold $M$ has a closed transverse  $(p,p)-$form $\Omega = \sigma_p \sum \Psi_A \wedge \overline{\Psi_A}$, where $\Psi_A$ is a holomorphic $p-$form (thus $\Omega \in (WP^p)^{int} \cap P^p$), then $M$ is called a {\bf holomorphically $p-$K\"ahler (pHK) manifold}.
\end{defn}

{\bf 6.2 Claim.} Every compact holomorphically parallelizable manifold of dimension $n$ is $(n-1)$HK, and it is K\"ahler if and only if it is a torus.

In fact $\Omega = \sigma_{n-1} \sum_{|A| = n-1} \varphi_A \wedge \overline{\varphi_A}$ is a closed strictly positive $(n-1,n-1)-$form; the second assertion is well known (see f.i. \cite{N}).
\medskip

We shall prove in the following theorem  (see f.i.  Theorem 3.2 in \cite{AB1}) that for this class of manifolds, all  type of K\"ahler properties coincide. This allows to use always {\it closed} currents, and moreover also {\it smooth} currents.

\begin{thm} The following properties are equivalent on a compact holomorphically parallelizable $n-$dimensional manifold, for every $p$, $1 \leq p \leq n-1$:

\begin{enumerate}
\item pHK
\item pK
\item pWK
\item pS
\item pPL
\item there are no strongly  positive currents $T \neq 0$, of bidimension $(p,p)$, such that $T = \ddb S$ for some {\bf smooth} current $S$ of bidimension $(p+1,p+1)$ 
\item there are no  simple exact holomorphic $(n-p)-$forms $\alpha \neq 0$. 
\end{enumerate}
\end{thm}

{\it Proof.} The only thing to prove is that (7) implies (1), because the other implications are true on every compact manifold (in particular, (6) implies (7) by Proposition 3.3). 

Let $\{\varphi_1, \dots, \varphi_n\}$ be a basis of $\Omega^1(M)$, and let 
$\Theta = \varphi_1 \wedge \dots \wedge \varphi_n$; consider the non-degenerate bilinear forms $F_j$ (to which we refer for orthogonality) given by
$$F_j : \Omega^j(M) \times \Omega^{n-j}(M)  \to \C, \ \  F_j (\beta, \rho)\Theta := \beta \wedge \rho.$$
Let $k=n-p$;  consider a basis $\{\rho_1, \dots, \rho_N\}$ of $\Omega^{k-1}(M)$, and let $\{d\rho_{s+1}, \dots, d\rho_N\}$ be the basis of the image of the operator $d: \Omega^{k-1}(M) \to \Omega^k(M)$. 

It is easy to check that all $\Psi \in (Im d)^{\perp}$ are closed. 
Thus $\Omega = \sigma_p \sum \Psi_h \wedge \overline{\Psi_h}$, where $\{\Psi_h\}$ is a basis for $(Im d)^{\perp}$, is the required pHK form. Indeed, it is closed and positive; it is transverse by Definition 2.4, because if 
$\Omega \wedge \sigma_k \varphi_I \wedge \overline{\varphi_I} = 0,$ for some $I$, $|I|=k$, then every $\Psi_h \wedge \varphi_I$ has to be zero, but  if $\varphi_I \in \Omega^k(M)$ is orthogonal to all $\Psi_h$, then $\varphi_I \in Im d$ would be a simple exact holomorphic $(n-p)-$form. 
\medskip

{\bf 6.4 Claim.} When $p < n-1$, and $M$ is a compact holomophically parallelizable pK manifold with a strictly {\it positive} pK form, then $M$ is K\"ahler (\cite{AB1}, Theorem 2.4). This result tells us that we can always choose a pK form in $(WP^p)^{int} \cap P^p$, but not  in $(P^p)^{int}$.
\medskip

 {\bf 6.5 $\eta \beta_{2n+1}$.} Among compact holomorphically parallelizable manifolds, we select suitable nilmanifolds which generalize the structure of $I_3$, the Iwasawa manifold, that is,  the manifolds $\eta \beta_{2n+1}$ of dimension $2n+1$  (see \cite{AB1}, Section 4). They are not \pkk for $p \leq n$ and are \pkk for $p > n$, and their forms  are in $(WP^p)^{int} \cap P^p$. 
 Let us describe them briefly.
 
 Let $G$ be the following subgroup of $GL(n+2, \C)$:
$$G := \{ A \in GL(n+2, \C) / A=\left(\begin{array}{ccc} 1 & X & z \\ 0 & I_n & Y \\ 0 & 0 & 1 \end{array} \right)
,  z \in \C, X, Y \in \C^n \},$$
and let $\Gamma$ be the subgroup of $G$ given by matrices with entries in $\Z[i]$. $\Gamma$ is a discrete subgroup and the homogeneous manifold 
$\eta \beta _{2n+1} := G/ \Gamma$ becomes a holomorphically parallelizable compact connected complex nilmanifold of dimension $2n+1$ (for $n=1$, $\eta \beta _{3}$ is nothing but the Iwasawa manifold $I_3$). 

The standard basis for holomorphic 1-forms on $\eta \beta _{2n+1}$ is called $\{ \varphi_1, \dots, \varphi_{2n+1} \}$, where the  $ \varphi_j$ are all closed, except $d \varphi_{2n+1}  =  \varphi_{1}  \wedge \varphi_{2}  + \dots + \varphi_{2n-1}  \wedge \varphi_{2n} .$
\medskip

\bigskip

{\bf \lq\lq Closure\rq\rq.} Let us consider now  the other condition used in the study of the product of \pkk manifolds: we always ask that the manifolds are \lq\lq $s-$K\"ahler\rq\rq  from a certain index $p$ on. This is not always true also starting from $p=1$, as we discuss now.

\medskip

Let $X$ be a complex manifold of dimension $n \geq 3$.
Consider the following questions, for $p > 1$:

When $\omega$ is a \lq\lq closed\rq\rq  strictly positive $(1,1)-$form, is $\Omega = \omega^p$  \lq\lq closed\rq\rq too?

 When $\Omega = \omega^p$ is \lq\lq closed\rq\rq, where $\omega$ is a  strictly positive $(1,1)-$form, is $\omega$ \lq\lq closed\rq\rq too?

\medskip

The answer to the first question is positive for the case pK, since
$ \de \omega^p = p \omega^{p-1} \wedge \de \omega.$
 
 This simply computation does not work in the pWK case and in the pPL case  (also for $p=n-1$), since in the first case 
 
 $\de \omega^p = p \omega^{p-1} \wedge \de \omega = p \omega^{p-1} \wedge \ddb \alpha$, whereas we need  $\de \omega^p = \ddb \beta$; 
 
 moreover,
 
 $\ddb \omega^p = p \omega^{p-2} \wedge ((p-1) \de \omega \wedge \overline{ \de \omega} + \omega \wedge \ddb \omega)$, and in particular 
 $\ddb \omega^2 = 2 ( \de \omega \wedge \overline{ \de \omega} + \omega \wedge \ddb \omega)$.
 
 Thus $\ddb \omega = 0$ does not imply $\ddb \omega^p =0$; but when both $\ddb \omega =0$ and $\ddb \omega^2 =0$, then $ \de \omega \wedge \overline{ \de \omega} =0$, so that $\ddb \omega^p =0$.
 \medskip

 In the pS case, if we consider the conditions given in Remark 3.1.1, we are in the same troubles as before. But when the question is translated as follows: \lq\lq Is a 1S manifold also a pS manifold?\rq\rq, the answer is positive.
 
 Indeed, let $\psi$ be a closed 2-form, whose $(1,1)-$component is $\omega > 0$. Then $\psi^p$ is closed too, and its $(p,p)-$component is given by $\omega^p + \zeta$, where $\zeta$ is a sum of $(p,p)-$forms of this kind: $\omega^j \wedge (\sigma_{p-j} \eta \wedge \overline{\eta})$, $\eta \in \Lambda^{p-j,0}$. Hence by Proposition 2.8, 
 $\zeta \in P^p$ is a positive form and $\omega^p$ is transverse, so that $\omega^p + \zeta$ is transverse.
 In fact, 
for every $V \in \Lambda_{p,0}(E)$, $V \neq 0$ and simple, we get 
$$(\omega^p + \zeta)(\sigma_p^{-1} V \wedge \overline V) > 
\zeta (\sigma_p^{-1} V \wedge \overline V) \geq 0.$$

 \medskip

As regards the second question, the answer is yes for the case pK, if $p < n-1$, since in this case
$d \omega^p =0 \iff d \omega =0$  by a direct computation. If $p=n-1$, it is false, indeed we have balanced non-K\"ahler manifolds.

In the other cases this does not work: for instance, $\ddb \omega^2 = 0$ does not imply $\ddb \omega = 0$, as some examples in \cite{FT2}, (2.3) show.  
 \bigskip

 {\bf 6.6 Construction.} Nevertheless, all examples known up to now are \lq\lq $s-$K\"ahler\rq\rq  from a certain index $p$ on. We give here a new example, a suitable nilmanifold  $M = G/\Gamma$ of complex dimension 4, to stress that this is not even the case.

We refer to \cite{EFV}, p. 216: in this paper the authors 
classify eight-dimensional nilpotent Lie algebras admitting an SKT (= 1PL)  structure. More precisely, in Theorem 4.2 they prove that, given a nilmanifold (different from a torus) $M = G/\Gamma$ with $dim_{\R} M = 8$, with an invariant complex structure $J$, there exists an SKT metric on $M$ compatible with $J$ if and only if the Lie algebra of $G$ belongs to one of the two families described in (4.1) and (4.3) ibidem.

As a matter of fact, while for the second family, as in the six-dimensional case, the SKT condition depends only on the complex structure (that is, given the invariant complex structure $J$ on $M = G/\Gamma$, the SKT condition is satisfied by either all invariant hermitian metrics, or by none), this is no more true in the first family (see Remark 4.1 ibidem).

Let us choose a very simple nilmanifold $(M,J)$ belonging to the first family of \cite{EFV}: that is, pick in (4.1):
$$B_1 = F_1 = B_4 = C_4 = 1, \ \ B_5 = C_3 = F_4 = F_5 = G_3 = G_4 = 0.$$
The constants satisfy (4.2), hence we get a 1PL metric on $M$. Let us give it explicitly:
\smallskip

Consider $M = G / \Gamma$, where $\Gamma$ is a lattice in the simply connected nilpotent Lie group $G$, and consider the
basis $\{\varphi_1, \varphi_2, \varphi_3, \varphi_4\}$ of invariant $(1,0)-$forms satisfying the following equations:
$$d \varphi_1 = d \varphi_2 = 0, \  \ \de \varphi_3 =  \varphi_1 \wedge \varphi_2, \ \  \db\varphi_3 = \varphi_1 \wedge \overline \varphi_1 +  \varphi_2 \wedge \overline \varphi_2, \ \ \de \varphi_4 =  \varphi_1 \wedge \varphi_2, \ \ \db\varphi_4 =0.$$
(Recall that if the structure equations are rational, then there exists a $\Gamma$ such that the quotient $M$ is compact, see \cite{EFV} p. 205).
Let 
$$\omega = \sigma_1 (\varphi_1 \wedge \overline \varphi_1 + \varphi_2 \wedge \overline \varphi_2 + \varphi_3 \wedge \overline \varphi_3 + \varphi_4 \wedge \overline \varphi_4)$$
 be the standard $(1,1)-$form: we get easily
$$i \ddb \omega =  i \ddb \sigma_1 (\varphi_3 \wedge \overline \varphi_3 + \varphi_4 \wedge \overline \varphi_4) = 0,$$
 thus $\omega$ is a 1PL form.

To check that $M$ is not 2PL, by Theorem 3.2 we need only to exhibit a strongly positive $(2,2)-$current $T \neq 0$ on $M$, with  $T = i \ddb A$ for a suitable $(1,1)-$current $A$. 
Let 
$$ A = -  i \ \varphi_3 \wedge \overline \varphi_3 ;$$
then 
$$i \ddb A =  \varphi_1 \wedge \varphi_2 \wedge \overline \varphi_1 \wedge \overline \varphi_2 := T \geq 0.$$

\medskip

{\bf 6.7 Example.} As regards the manifold $M$, the situation is the following: 
\begin{enumerate}

\item $M$ is 1PL due to the $(1,1)-$form $\omega$;

\item $M$ is not 2PL, as shown above; thus it is not  2S, 2WK and 2K;

\item $M$ is not 1K nor 1S, since this would contradict (2). Hence it is not 2WK;

\item $M$ is 3PL by Proposition 3.3, but it is not 3S (hence not 3WK nor 3K) since
$$T= dR = d(i(\varphi_3 - \varphi_4) -i(\overline \varphi_3 - \overline \varphi_4)) = 2i(\varphi_1 \wedge \overline \varphi_1 +  \varphi_2 \wedge \overline \varphi_2) \geq 0.$$
\end{enumerate}

\bigskip

To end our paper let us give here a couple of examples, to illustrate the results we got on the product of \pkk manifolds.
\medskip

 In \cite{AB} we build an example to show that, when $f:\tilde X \to X$ is a proper  modification, we can sometimes pull-back  \pkk properties for $p>1$. Indeed, we consider  a smooth modification $\tilde X$ of $\P_5$, where the center $Y$ is a surface with a singularity; the singular fibre has two irreducible components, one of which is biholomorphic to $\P_2$ and the other is a holomorphic fibre bundle over $\P_1$ with $\P_2$ as fibre. We show that $\tilde X$ is not K\"ahler, because it contains a copy of Hironaka's threefold, but it is \pkk for every $p > 1$; moreover, its $pK-$forms are strictly strongly positive (ibidem, p. 313).
\medskip

{\bf 6.8 Example.} Let us fix $m \geq 1$, and consider $(\P_1)^m \times \tilde X$.

By Theorem 5.5, $(\P_1)^m \times \tilde X$ is sK for $s \geq m+2$; it is not K\"ahler because it contains $\tilde X$ and it is not pK, for $2 \leq s \leq m+1$ by Corollary 4.4.
\bigskip

To give another interesting example, let us consider the  deformations of nilmanifolds (see f.i. \cite{N}), in particular of the Iwasawa manifold $I_3$. Recall that we can find on $I_3$ a basis $\{\varphi_1, \varphi_2, \varphi_3\}$ of holomorphic 1-forms, such that
$$d \varphi_1 = d \varphi_2 = 0, \ d\varphi_3 =\de \varphi_3 =  \varphi_1 \wedge \varphi_2.$$

In \cite{AB2} we studied a small deformation of $I_3$, namely $I_{3,t}$ with $|t|$ small, which has a basis of invariant $(1,0)-$forms such that 
$$d \varphi_1 = d \varphi_2 = 0, \  \ \de \varphi_3 =  \varphi_1 \wedge \varphi_2, \ \  \db\varphi_3 = - t \varphi_2 \wedge \overline \varphi_2.$$
When $t \neq 0$, $I_{3,t}$ is no more 2K.
\medskip

In \cite{FPS} the authors study SKT metrics on (real) six-dimensional nilmanifolds, and prove that the condition $\ddb \omega =0$ only depends on the underlying complex structure. In particular, they give in Theorem 1.2 a necessary and sufficient condition to the existence of a SKT metric in terms of the differentials of a basis of invariant $(1,0)-$forms $\{\varphi_1, \varphi_2, \varphi_3\}$. This holds in particular for 
$I_3$; choosing (in their notation) $A=D=0, E=-1, C = i/2, B = -i$, that is $t=s=1$, we get a manifold
$I_{3,1}$ which has a basis of invariant $(1,0)-$forms such that 
$$d \varphi_1 = d \varphi_2 = 0, \  \ \de \varphi_3 =  \varphi_1 \wedge \varphi_2, \ \  \db\varphi_3 = \frac{i}{2} (\varphi_1 \wedge \overline \varphi_1 + 2 \varphi_2 \wedge \overline \varphi_2).$$

Straightforward computations show that indeed the transverse $(1,1)-$form $\omega = \sigma_1 (\varphi_1 \wedge \overline \varphi_1 +  \varphi_2 \wedge \overline \varphi_2 + \varphi_3 \wedge \overline \varphi_3)$ is $\ddb-$closed.
\medskip

We shall now compute all K\"ahler degrees of $I_{3,1}$. For simplicity, denote $\varphi_{j,k} := \varphi_j \wedge \varphi_k, \ \overline \varphi_{j,k} := \overline \varphi_j \wedge \overline{\varphi_k}, $ 
$\varphi_{j,k,h} := \varphi_j \wedge \varphi_k \wedge \varphi_h$ and so on.
  
First of all, consider $\omega^2 = 2 \sigma_2 (\varphi_{1,2} \overline{\varphi_{1,2}} + \varphi_{1,3} \overline{\varphi_{1,3}} + \varphi_{2,3} \overline{\varphi_{2,3}})$; we get 

$$\de \omega^2 = \frac{-3i}{4} \varphi_{1,2,3} \overline{\varphi_{1,2}} = - \db (\frac{3i}{4} \varphi_{1,2,3} \overline{\varphi_{3}})$$
with $\de (\frac{3i}{4} \varphi_{1,2,3} \overline{\varphi_{3}}) = 0$. By Remark 3.1.1, this implies that $I_{3,1}$ is 2S (and obviously 2PL).

Let
$$\de \overline \varphi_{3} + \db \varphi_{3} = i (\varphi_1 \wedge \overline \varphi_1 + 2 \varphi_2 \wedge \overline \varphi_2) := T \geq 0,$$
 where $\ddb \varphi_{3} =0$: therefore by Theorem 3.2 $I_{3,1}$ is not 2WK (nor 2K). Notice that it is known that the same metric cannot be SKT and balanced (see f.i. \cite{FPS}, Proposition 1.4), but here we get more: $I_{3,1}$ cannot have balanced metrics.
 
In the same manner, consider the $(2,1)-$current given by $R := \varphi_{1,2} \overline{\varphi_{3}}$; since $\de R =0$, we get 
$$\de \overline R + \db R = d(R + \overline R) =  2 \varphi_{1,2} \overline{\varphi_{1,2}} := T \geq 0,$$
so that $I_{3,1}$ is not 1S (nor 1K and 1WK).
\smallskip
 
 With the same kind of computations, one can prove that $I_{3,t}$ is not \kk, and is not 2WK (nor 2K), but is 2S and 2PL.
\medskip

{\bf 6.9 Example.} Consider $X := I_{3,1}, Y := \eta \beta_5$; then $M = I_{3,1} \times \eta \beta_5$ is
\begin{enumerate}
\item 7PL by Proposition 3.3;

\item 6PL by Theorem 5.4;

\item 7S by Proposition 5.2;

\item not \kk nor \lq\lq 2-K\"ahler\rq\rq by Corollary 4.2;

\item not \lq\lq 4-K\"ahler\rq\rq nor \lq\lq 5-K\"ahler\rq\rq by Corollary 4.4;

\item not 6K, 6WK, 6S, 7K, 7WK  by Corollary 4.4;

\item not \lq\lq 3-K\"ahler\rq\rq by Proposition 3.3; in fact, pull-back on $I_{3,1} \times \eta \beta_5$ the forms $\varphi_j$ of $I_{3,1}$ and $\varphi_k'$ of $\eta \beta_5$  and notice that $\varphi_1 \wedge \varphi_2 \wedge \varphi_1' \wedge \varphi_3' \wedge \varphi_4'  = \de (\varphi_1 \wedge \varphi_2 \wedge \varphi_1' \wedge \varphi_5')$ 
is a simple exact holomorphic 5-form.
\end{enumerate} 

\bigskip
\bigskip

\end{document}